\documentclass[10pt]{amsart}
\newtheorem{thm}{Theorem}[section]
\newtheorem{prop}[thm]{Proposition}
\newtheorem{cor}[thm]{Corollary}
\newtheorem{lem}[thm]{Lemma}
\newtheorem{dfn}[thm]{Definition}
\newtheorem{remark}{Remark} 
\def\bi{{\bf i}}
\def\b1{{\bf 1}}
\def\bR{{\bf R}}
\def\bC{{\bf C}}
\def\bZ{{\bf Z}}

\begin{document}
{
\center{\Large{\bf
Invariant of the hypergeometric group
associated to the quantum cohomology of the projective space.
} }
 \vspace{1pc}
\center{\large{ Susumu TANAB\'E}}
 }
 \vspace{0.5pc}

\begin{center}
\begin{minipage}[t]{12.2cm}
{\sc Abstract -} {\em
We present a simple method to calculate the Stokes matrix for the
quantum cohomology of the projective spaces ${\bf CP}^{k-1}$
in terms of certain hypergeometric group.
We present also an algebraic variety whose
fibre integrals are solutions
to the given hypergeometric equation.
 }
\end{minipage}
\end{center}

{
\center{\section{\bf Generalized hypergeometric function
}}
}
We begin with a short review on  the motivation of our problem making
reference to the works ~\cite{Dub}, ~\cite{Guz} where one can find
precise definitions of the notions below.

At first, we consider a $k-$ dimensional Frobenius manifold $F$ with flat
coordinates
$(t_1,\cdots, t_k)$$\in F$ where the coordinate $t_i$
corresponds to coefficients of the basis $\Delta_i$
of the quantum cohomology $H^{\ast}({\bf CP}^{k-1}).$
On $H^{\ast}({\bf CP}^{k-1})$
one can define so called quantum multiplication
$$ \Delta_\alpha \bullet \Delta_\beta =
C^{\gamma}_{\alpha,\beta}\Delta_\gamma,$$
or $$\frac{\partial} {\partial t_\alpha} \cdot \frac{\partial} {\partial
t_\beta}=  C^{\gamma}_{\alpha,\beta}\frac{\partial}
{\partial t_\gamma},$$
on the level of vector fields on $F.$ The Frobenius manifold is furnished with the Frobenius algebra on the tangent space $T_t F$ depending analytically on
$t \in F,$ $T_t F= (A_t, <\;,\; >_t)$ where $A_t$ is a commutative associative
$\bf C$ algebra and $ <\;,\; >_t:A_t \times A_t \rightarrow \bf C$
a symmetric non-degenerate bilinear form. The bilinear form  $ <\;,\; >_t$
defines a metric on $F$ and the Levi-Civita connexion $\nabla$ for this metric can be considered. Dubrovin introduces a deformed flat connexion
$\tilde \nabla$ on $F$ by the formula $\tilde \nabla_u v :=\nabla_u v+
xu \cdot v $ with $x \in \bf C$ the deformation parameter.
Further he extends $\tilde \nabla$ to $F \times \bf C$. Especially we have
$\tilde \nabla_{\frac{\partial} {\partial x} } = \frac{\partial} {\partial x}
- E(t) - \frac{\mu}{x},$
where $E(t)$ corresponds to the multiplication by the Euler vector field
$ E(t) = \sum_{1 \leq j \not = 2 \leq k-1}(2-j)
t_j\frac{\partial} {\partial t_j} + kt_2\frac{\partial} {\partial t_2}.$

After ~\cite{Dub},~\cite{Guz} the quantum cohomology $\vec{u}(x)
=(u_1(x),\cdots, u_k(x))$
for the projective space ${\bf CP}^{k-1}$  at a semisimple point
$(0,t_2,0,\dots,0)$ (i.e. the algebra $(A_t, <\;,\; >_t)$ is semisimple there)
satisfies the following system of differential equation:
$$\partial_x \vec{u}(x)=
(k \mathcal{C}_2(t)+\frac{\mu}{x})\vec{u}(x).\leqno(1.1)$$
where
$$\mathcal{C}_2(0,t_2,0,\dots,0)=
(C^{\gamma}_{2,\beta})_{1 \leq \beta,\gamma\leq k}=
\left(
\begin{array}{llcll}
0 & 0 & \cdots &0 &e^{t_2} \\
1 & 0 &  \ddots &0 &0 \\
0 &1 & \ddots & 0&0 \\
\vdots &\ddots  & \ddots &\vdots &\vdots\\
0 & 0 & \cdots &  1 &0 \\
\end{array} \right).$$
The matrix $\mu$ denotes a diagonal matrix
with rational entries:
$$ \mu = diag \{-\frac{k-1}{2}, -\frac{k-3}{2},\cdots,
\frac{k-3}{2},\frac{k-1}{2} \}.$$

The last  component $u^k(z)$ (after the change of variables $z := kx
e^{\frac{t_2}{k}}$)
of the above system for the quantum cohomology
satisfies a differential equation
as follows ~\cite{Guz}:
$$[(\vartheta_z)^k- z^k]z^{\frac{-k+1}{2}} u^k(z)=0, \leqno(1.2)$$
with $\vartheta_z =z \frac{\partial} {\partial z}.$
After the Fourier-Laplace transformation
$$\tilde u(\lambda) = \int e^{\lambda z}z^{\frac{-k+1}{2}}u^k(z) dz,$$
we obtain an equation as follows:
$$[(\vartheta_ \lambda+1)^k-(\frac{\partial}
{\partial \lambda})^k]\tilde u(\lambda) =0.$$
Here the notation $\vartheta_ \lambda$ stands for
$\lambda\frac{\partial} {\partial \lambda}.$
After multiplying $\lambda^k$ from the left, we obtain
$$[\lambda^k(\vartheta_ \lambda+1)^k-
\vartheta_\lambda(\vartheta_\lambda-1)
(\vartheta_\lambda-2)\cdots (\vartheta_\lambda-(k-1))]
\tilde u(\lambda) =0.$$
The equation for $\lambda \tilde u(\lambda ),$ the Fourier-Laplace
transform of $\frac {\partial}{\partial z} z^{\frac{-k+1}{2}}u^k(z)$ should be
$$ [\lambda^k(\vartheta_ \lambda)^k-
(\vartheta_\lambda-1)
(\vartheta_\lambda-2)\cdots (\vartheta_\lambda-k)]
(\lambda \tilde
u(\lambda)) =0.\leqno(1.3)$$
It is evident that the Stokes matrix for $\frac {\partial}{\partial z}
z^{\frac{-1+k}{2}}u^k(z)$ is identical with that of the original solution $u^k(z).$

Before proceding further, we remind the following theorem that gives connexion
between the Stokes matrix of the system $(1.1)$
with the monodromy of the equation
$(1.3).$
Let us consider the Fourier-Laplace transform of the system $(1.1)$:
$$(\vartheta_\lambda +id_k)\vec{\tilde u}(\lambda)=
(k \mathcal{C}_2(t)\partial_\lambda - \mu)
\vec{\tilde u}(\lambda).\leqno(1.1)'$$

In a slightly more general setting, let us observe a system with
regular singularities:
$$
(\Lambda - \lambda\cdot id_k)\partial_\lambda
\vec{\tilde u}(\lambda)=(id_k + A_1(\lambda))\vec{\tilde u}(\lambda)
\leqno(1.1)"$$
with $\Lambda \in GL(k, {\bf C})$ whose eigenvalues $(\lambda_1, \cdots,
\lambda_{k})$ are all distinct,
$A_1(\lambda) \in End({\bf C}^k) \otimes {\mathcal O}_{\bf C}$ with
$A_1(0)= diag (\rho_1, \cdots, \rho_k)$
where none of the $\rho_j$'s is an integer.
We call solutions to a scalar differential equation deduced from
$(1.1)"$ {\it component solutions}. Thus solutions to $(1.3)$ are
component solutions to $(1.1)".$

 \begin{thm} (~\cite{BJL}, ~\cite{Dub}) Under the assumption
that the eigenvalues of the matrix $A_1(0)$ are distinct, the Stokes
matrix $S$ for the component solutions of $(1.1)$ expresses the
symmetric Gram matrix $G$ of the component solutions of $(1.1)'$ as
follows:  $$ ^tS+S=2G.$$ \label{thm11} \end{thm}

As for the definition
of the Stokes matrix $S$ for the system $(1.4)$ we refer to
~\cite{Dub}, ~\cite{Guz}.  The main theorem of this article is the
following:  \begin {thm} The $i,j$ component $S_{ij}$, $1 \leq i,j\leq
k$ of the Stokes matrix to the system $(1.1)$ has the following
expression:  $$  S= \left\{ \begin{array}{ll}
(-1)^{i-j}\;_kC_{i-j}\;\;\;\;& i \geq j\\
0\;\;\;\; &i < j\\
\end{array}
\right.\;
$$
\label{thm12}
\end{thm}
This theorem has already been shown by D.Guzzetti ~\cite{Guz}
by means of a detailed study of braid group actions etc on the set of
solutions to $(1.2).$ We present here another approach to understand
the structure of the Stokes matrix.

\begin{remark}
In this article we observe the convention of the matrix multiplicatin as follows:
$$ A\cdot x= (a_{ij})_{0 \leq i,j \leq k-1} (x_i)_{0 \leq i\leq k-1} = \langle \sum_{i=0}^{k-1}
a_{ij} x_i \rangle_{0 \leq j \leq k-1}.$$
The matrix operates on the vector from left, in contrast to the convetion used in ~\cite{Dub},~\cite{Guz}.
\label{remark1}
\end{remark}
On the other hand it has been known since ~\cite{Beil} that a collection
of coherent sheaves ${\mathcal O}(-i)$
$0 \leq i \leq k-1$ on ${\bf CP}^{k-1}$
satisfies the following relation
$$ Hom({\mathcal O}(-i), {\mathcal O}(-j)) = S^{i-j}({\bf C}^k), 0 \leq i,j
\leq k-1$$
$$Ext^\ell({\mathcal O}(-i),{\mathcal O}(-j))=0, 0 \leq i,j
\leq k-1, \ell >0$$
These relation entails immediately the equality
$$\chi( {\mathcal O}(-i), {\mathcal O}(-j)):=\sum_{\ell=0} (-1)^{\ell}
Ext^\ell({\mathcal O}(-i),{\mathcal O}(-j)) =
\left\{
\begin{array}{ll}
\;_{k+i-j-1}C_{i-j}\;\;\;\;& i \geq j\\
0\;\;\;\; &i < j.\\
\end{array}
\right. $$
We consider action of the braid group $\beta_i  \in {\bf B}_k, 1 \leq i
\leq k-1$ that corresponds to the braid action between $i-$th basis and
$(i+1)-$st basis of the space on which act a matrix.
In our situation, $\beta_i$ represents the braid action between
${\mathcal O}(1-i)$ and ${\mathcal O}(-i).$ In literature on
coherent sheaves on algebraic varieties, this procedure is called
mutation (e.g. ~\cite{GorKul}).
Let us denote by $\beta$ an element of the braid group ${\bf B}_k$
$$ \beta = \beta_1(\beta_2 \beta_1) \cdots (\beta_{k-1} \cdots \beta_2
\beta_1).$$ We introduce a matrix of reordering $J = \delta_{i,
k-1-i}, 0 \leq i \leq k-1.$ In this situation our Stokes matrix from
Theorem~\ref{thm12} is connected with the matrix ${\bf \chi} := (\chi(
{\mathcal O}(-i), {\mathcal O}(-j))$ $0 \leq i,j \leq k-1$  in the
following way, $$ ^tS= J \beta \chi \beta J.$$ Eventually it turns out
that $ \chi = S^{-1}.$ This general fact on the braid group is
explained in ~\cite{Z}, \S 2.4.

As our Stokes matrix is determined up to the change of basis,
including effects by braid group actions,
the Theorem ~\ref{thm12} is a confirmation of an hypothesis \cite{Dub2}
that the matrix for certain exceptional collection of coherent sheaves
on a good Fano variety $Y$ must coincide with the Stokes matrix for the
quantum cohomology of $Y.$

Our  strategy to prove Theorem ~\ref{thm12} consists in the study of
system $(1.1)'$, instead of $(1.1)$ itself.

Further we consider so called the Kummer covering (naming after N.Katz)
of the projective space ${\bf CP}^1$ by $\zeta = \lambda^k$ to
deduce an hypergeometric equation:
$$[\zeta(\vartheta_ \zeta)^k-
(\vartheta_\zeta-\frac{1}{k})
(\vartheta_\zeta-\frac{2}{k})\cdots (\vartheta_\zeta-1)]
v(\zeta) =0, \leqno(1.4)$$
for $v(\lambda^k)= \lambda \tilde u(\lambda).$
We remind of us  here a famous theorem due to A.H.M.Levelt
 that allows us to express the (global) monodromy group
of the solution to $(1.4)$ in quite a simple way.
For the hypergeometric equation in general,
$$ [\prod_{\ell=1}^k (\vartheta_\zeta-\alpha_\ell) -
\zeta\prod_{\ell=1}^k (\vartheta_\zeta-\beta_\ell)]v(\zeta)=0, \leqno(1.5)$$
we define two vectors $(A_1, \cdots, A_k)$ and $(B_1, \cdots, B_k)$
in the following way:
$$\prod_{\ell=1}^k (t-e^{2\pi \alpha_\ell i})
=t^k+ A_1t^{k-1}+A_2t^{k-2}+ \cdots + A_k,$$
$$\prod_{\ell=1}^k (t-e^{2\pi \beta_\ell i})
=t^k+ B_1t^{k-1}+B_2t^{k-2}+ \cdots + B_k.$$
\begin{dfn}
A linear map $L \in GL(k, {\bf C})$
is called pseudo-reflexion if it satisfies the condition
$rank(id_k- L)=1.$ A pseudo-reflexion $R$ satisfiying an additional condition
$R^2=id_k$ is called a reflexion.
\label{dfn11}
\end{dfn}

\begin{prop} (~\cite{BH},~\cite{Lev}) For the solutions to $(1.5)$, the
monodromy action on them at the points $\zeta =0,\infty, 1$
has the follwing expressions:
$$ h_0=
\left(
\begin{array}{llccll}
0 & 0 & \cdots &0 & -A_k \\
1 & 0 &  \ddots &0 &-A_{k-1} \\
0 &1 & \ddots  &0&-A_{k-2} \\
\vdots &\ddots  & \ddots &\vdots &\vdots\\
0 & 0 & \cdots &  1 &-A_1 \\
\end{array} \right),
\leqno(1.6)$$
$$ (h_\infty)^{-1}=
\left(
\begin{array}{llccll}
0 & 0 & \cdots &0 & -B_k \\
1 & 0 &  \ddots &0 &-B_{k-1} \\
0 &1 & \ddots  &0&-B_{k-2} \\
\vdots &\ddots  & \ddots &\vdots &\vdots\\
0 & 0 & \cdots &  1 &-B_1 \\
\end{array} \right),
$$
whereas $h_1 = (h_0 h_\infty)^{-1}$ is a pseudo-reflexion.
\label{prop12}
\end{prop}

It is worthy to notice that the above
proposition does not precise for which
bases of solution to $(1.5)$ the monodromy is calculated.
  As a corollary to the Proposition ~\ref{prop12}, however,  we
see that the monodromy action on the solutions to our equation
$(1.4)$ can be written down with respect to a certain basis as follows:
$$h_0=\left(
\begin{array}{llcll}
0 & 0 & \cdots &0 &1\\
1 & 0 &  \ddots &0 &0 \\
0 &1 & \ddots & 0&0 \\
\vdots &\ddots  & \ddots &\vdots &\vdots\\
0 & 0 & \cdots &  1 &0 \\
\end{array} \right),\leqno(1.7)$$
$$h_\infty =\left(
\begin{array}{llcclll}
\;_kC_1 &1&0&\cdots  & 0&0  \\
\;-_{k}C_2 & 0&1 &\ddots & 0 &0 \\
\vdots & \vdots&\vdots & \ddots &\vdots &\vdots \\
(-1)^{k-1}\;_kC_{k-2} & 0 &0& \ddots &1&0\\
(-1)^{k}\;_kC_{k-1} & 0 & 0&\cdots &  0 &1 \\
-(-1)^k&0 &0& \cdots &0 &0
\end{array} \right).$$
In other words,
$$ det(t-h_0) = t^k-1, det(t-h_\infty) = (t-1)^k. \leqno(1.8)$$
Furthermore we have,
$$h_1 =\left(
\begin{array}{llcclll}
\;(-1)^{k-1} &0&0&\cdots  & 0&0  \\
 (-1)^{k-2}\;_kC_{k-1} & 1&0 &\ddots & 0 &0 \\
(-1)^{k-3}\;_kC_{k-2} & 0 &1& \ddots &0&0\\
\vdots & \vdots&\vdots & \ddots &\vdots &\vdots \\
\;_kC_{1}&0 &0& \cdots &0 &1
\end{array} \right). \leqno(1.9)$$
In the next sectin we will see that the theory of Levelt supplies us with
necessary data to calculate further the Stokes matrix of the
solutions to $(1.1).$

{
\center{\section{\bf
Invariants of the Hypergeometric group
}}
}
Let us begin with a detailed description of the generators of the
hypergeometric group defined for the solutions to the  equation (1.3).

\begin{prop}(cf. \cite{Gol} I, 8.5)
The generators of the hypergeometric group $H$  of the equation(1.3)
are expressed in terms of the matrices introduced in the
Proposition~\ref{prop12} as follows:
$$ M_0= h_0^k=1, M_1 = h_1= (h_0 h_\infty)^{-1}, M_\infty = h_\infty^k,
M_{\omega^i} = h_\infty^{-i}h_1 h_\infty^i (i=1,2, \cdots, k-1), \leqno(2.1)
$$
where $M_t$ denotes the monodromy action around the point
$t\in {\bf CP}^1_{\lambda}.$  The generators around singular points
$\omega^i =e^{2 \pi \sqrt -1 \frac{i}{k}}$  naturally
satisfy the Riemann-Fuchs relation:  $$M_\infty M_{\omega^{k-1}}
M_{\omega^{k-2}}\cdots M_{\omega} M_1 = id_k. \leqno(2.2)$$
\label{prop13}
\end{prop}
\noindent
{\bf proof} Let us think of a $k-$leaf covering
$\tilde {\bf CP}^1_\lambda$ of ${\bf CP}^1_\zeta$
that corresponds to the Kummer covering $\zeta^k =\lambda.$
In lifting up the path around $\zeta=1$ the first leaf of
$\tilde {\bf CP}^1_\lambda,$ the monodromy $h_1$ is sent
to the conjugation with a path around $\lambda=\infty.$ That
is to say we have $M_{\omega}= h_\infty^{-1}h_1 h_\infty.$
For other leaves the argument is similar. {\bf Q.E.D.}

Let us denote by $X^K$ a $k\times k$ matrix that satisfies
the relation
$$ {\bar g}X^K\;^t{g}=X^K, \leqno(2.3)$$
for every element $g$ of a group $K\subset GL(k,{\bf C}).$
From the definition, the set of all $X^K$ for a group $K$
reperesents a $\bf C$ vector space in general. We will call
a matrix of  this space the quadratic invariant of the
group $K.$
In the special case in which we are interested,
the following statement holds.

\begin{lem}
For the hypergeometric group $H$ generated by the pseudo-reflexions as in
(2.1), for every $X^H$ there exists a non zero $k\times k $ matrix
${\tilde X}^H$ such that $X^H =\lambda {\tilde X}^H $ for some
$\lambda \in {\bf C}\setminus \{0\}.$
\label{lem22}
\end{lem}
\noindent
{\bf proof}
The relation $$h_1 X\;^t{h_1}=X \leqno(2.4)$$
gives rise to equations on
$x_{0j}$ and $x_{j0}.$ That is to say, the first row of
(2.4) corresponds to
$$ (-1)^i \; _kC_{i} x_{00} - (-1)^{k-1}x_{0i}=x_{0i}, \;\; 1 \leq i \leq k-1,
$$  while
$$ (-1)^i \; _kC_{i} x_{00} - (-1)^{k-1}x_{i0}=x_{i0}, \;\; 1 \leq i \leq k-1.
$$ Thus we obtained $2(k-1)$ linearly independent equations.
Further
by concrete calculus one can easily see that
$$M_{\omega^\ell} = id_k + T_\ell,$$
where
$$T_\ell = \left(
\begin{array}{lllclll}
t^{(\ell)}_{0}\tau_0 & t^{(\ell)}_{1}\tau_0 &\cdots  &t^{(\ell)}_{\ell}\tau_0&
0&\cdots &0\\
t^{(\ell)}_{0}\tau_1 & t^{(\ell)}_{1}\tau_1 &\cdots  &t^{(\ell)}_{\ell}\tau_1&
0&\cdots &0\\
\vdots &\vdots &\vdots & \vdots &\vdots &\vdots &\vdots\\
t^{(\ell)}_{0}\tau_{k-1} & t^{(\ell)}_{1}\tau_{k-1} &\cdots
&t^{(\ell)}_{\ell}\tau_{k-1}&
0&\cdots &0
\end{array} \right),$$
with $(k-\ell)- $ zero columns from the right. The remaining columns are
generated from
$T_1$ after simple linear recurrent relations by an inductive way.
The relation
 $M_{\omega} X\; ^tM_{\omega}=X$
gives rise to new equations
$$ (1+t^{(1)}_{1}\tau_1)^2 x_{11} + \rm{linear\; functions\; in}\;\it
(x_{0i}, x_{i0}  )=x_{11},$$
with $(1+t^{(1)}_{1}\tau_1) = -1 + (_kC_1)^2 \not = 1$
and
$$ (1+t^{(1)}_{1}\tau_1) x_{1i} + \rm{linear\; functions\; in}\;\it
(x_{0i}, x_{i0}, x_{11}  )=x_{1i},$$
$$ (1+t^{(1)}_{1}\tau_1) x_{i1} + \rm{linear\; functions\; in}\;\it
(x_{0i}, x_{i0}  )=x_{i1}.$$
Thus we get $2k-3$ new linearly independent equations.
In general for $(\ell, \ell)$ term, we get from the relation
$M_{\omega^\ell} X\; ^tM_{\omega^\ell}=X, 1 \leq \ell \leq k-1$
$$(1+t^{(\ell)}_{\ell}\tau_\ell)^2 x_{\ell \ell} + \rm{linear\; functions\; in}\;\it
(x_{\nu i}, x_{i \nu}, 0 \leq \nu \leq \ell-1)=x_{\ell \ell},$$
with $1+t^{(\ell)}_{\ell}\tau_\ell= -1 + (_kC_\ell)^2 \not = 1.$
For $x_{i\ell}$
$$ (1+t^{(\ell)}_{\ell}\tau_\ell) x_{i\ell} + \rm{linear\; functions\; in}\;\it
(x_{\nu i}, x_{i \nu}, 0 \leq \nu \leq \ell-1, x_{\ell \ell})=x_{i\ell}.$$
In this way we get a set of
$2(k-1)+\sum_{\ell=1}^{k-1}\bigl( 2(k-\ell)-1\bigr)=
k^2-1$
independent linear equations with respect to
the elements of $X.$ {\bf Q.E.D.}

The quadratic invariant $X^{H_0}$ for $H_0 =\{h_0, h_\infty\}$ is
invariant with respect to $H.$  After the lemma ~\ref{lem22},
$\bC$ vector space of quadratic invariants $X^H$ is one dimensional.
Thus every $X^{H_0}$ is also $X^{H}.$
Hence we can
calculate the quadratic invariant $X^H$ after the following relations,
$$\bar{h_0}X^H\;^th_0=X^H,
\bar{h_{\infty}}X^H\;^t{h_\infty}=X^H.\leqno(2.5)$$

    From \cite{BH} we know that the inverse to $X^{H_0}=X^{H}$, if it exists,  must be  a Toeplitz matrix i.e.:
$$(X^{H_0})^{-1}=\left(
\begin{array}{llcll}
x_{0} & x_{1} & x_{2} &\cdots &x_{k-1}\\
x_{-1} & x_{0} &  x_{1} & \cdots &x_{k-2} \\
x_{-2} &x_{-1} & x_{0} &  \cdots&x_{k-3} \\
\vdots &\ddots  & \ddots &\vdots &\vdots\\
x_{-(k-1)} & x_{-(k-2)}& x_{-(k-3)}& \cdots &x_{0} \\
\end{array} \right).$$
Making use of this circumstances, it is possible to show
that the system of equations that arises
from the relations
$$ ^th_{\infty} (X^{H_0})^{-1} \;{\bar h}_{\infty}=(X^{H_0})^{-1},
^th_{0} (X^{H_0})^{-1} \;{\bar h}_{0}=(X^{H_0})^{-1}
,$$
for $(X^{H_0})^{-1}$
consists of $2(k-1)$ equations.
$$x_{k-1-i}= x_{-i-1},\leqno(2.6)'$$
$$(-1)^{k+1}x_{k-1-i}+ (-1)^k \;_kC_{k-1}x_{k-2-i}+ \cdots
+\; _kC_3x_{2-i}-\;_kC_2x_{1-i}+ k x_{-i}= x_{-1-i} \leqno(2.6)''$$
This calculates  the matrix $X^H$ for the case $k-$ odd.

As for the case $k-$ even,
our matrix $X^H$ has the following form
$$X^{H}=\left(
\begin{array}{llcll}
0 & y_{1} & y_{2} &\cdots &y_{k-1}\\
y_{-1} & 0 &  y_{1} & \cdots &y_{k-2} \\
y_{-2} &y_{-1} & 0 &  \cdots&y_{k-3} \\
\vdots &\ddots  & \ddots &\vdots &\vdots\\
y_{-(k-1)} & y_{-(k-2)}& y_{-(k-3)}& \cdots &0 \\
\end{array} \right),$$
where $y_{-(k-1)}, \cdots y_{k-1}$ satisfy $2(k-1)$
equations for some constant $y_0,$
$$ y_i + y_{-i}=0, \;\; y_i - y_{-i}= \;
2(-1)^i\;_kC_i y_0\;\; 1 \leq i \leq k-1, \leqno(2.6)'''$$
which are derived from (2.5).
Thus the matrix $X^H$ for the case $k-$even is
obtained.

We remember here a classical theorem on the pseudo-reflexions.

\begin{thm}(cf. Bourbaki Groupe et Alg\`ebre de Lie Chapitre V,
\S 6, Exercise 3)
Let $E$ be a vector space with basis $(e_1, \cdots, e_d),$
and their dual basis $(f_1, \cdots, f_d)$$\in E^\ast.$
Let us set $a_{ij}= f_i(e_j).$ The pseudo-reflexion $s_i$
with respec to the basis $f_i$ is defined as
$$s_i(e_j)=  e_j - f_i(e_j) e_i= e_j - a_{ij} e_i.$$
Set $$V=\left(
\begin{array}{llcll}
a_{11} & a_{21} & a_{31} &\cdots &a_{d1}\\
0 & a_{22} &  a_{32} & \cdots &a_{d2} \\
0 &0 & a_{33} &  \cdots&a_{d3} \\
\vdots &\ddots  & \ddots &\vdots &\vdots\\
0 & 0 & \cdots &  0 &a_{dd} \\
\end{array} \right),
U=\left(
\begin{array}{llcll}
0 & 0 & 0 &\cdots &0\\
a_{12} & 0&  0& \cdots &0 \\
a_{13} &a_{23} & 0 &  \cdots&0 \\
\vdots &\ddots  & \ddots &\vdots &\vdots\\
a_{1d} & a_{2d} & \cdots & a_{d-1,d}&0 \\
\end{array} \right) \leqno(2.7)$$
Under these notations, the composition of all possible
reflexions $s_d s_{d-1} \cdots s_1$
(a Coxeter element) with
respect to the basis
$(e_1, \cdots, e_d)$
is expressed as follows:
$$s_d s_{d-1} \cdots s_1= (id_d - V)(id_d+U)^{-1}. \leqno(2.8)$$
\label{thm23}
\end{thm}
\noindent
{\bf proof}
For $1 \leq i,k \leq d$ we define
$$ y_i = s_{i-1} \cdots s_{1}(e_i).$$
It is possible to see that $$e_i= y_i+\sum_{k<i \leq d}a_{ki}y_k,
s_d \cdots s_1(e_i)=y_i - \sum_{i\leq k \leq d}a_{ki}y_k.$$
The statement follows immediately from these relations.
{\bf Q.E.D.}

To establish a relationship between the invariant $X^H$ and the Gram
matrix necessary for calculus of the Stokes matrix, we
investigate a situation
where the generators of the hypergeomeric group have special forms. Namely
consider a hypergeometric group $\Gamma$ of rank $k$ generated by
pseudo-reflexions
$R_0, \cdots R_{k-1}$ where
$$R_j= id_k -  Q_j, \leqno(2.9)$$
with
$$Q_j = \left(
\begin{array}{llcllll}
0 & \cdots &0& t_{j0}&0 &\cdots &0\\
0 & \cdots &0& t_{j1}&0 &\cdots &0\\
0 & \cdots &0& t_{j2}&0 &\cdots &0\\
\vdots &\cdots& \vdots &\vdots
&\vdots&\cdots &\vdots\\
0 & \cdots &0& t_{j,k-1}&0 &\cdots &0\\
\end{array} \right), \; 0 \leq j \leq k-1, \leqno(2.10)$$
all zero components except for the $j-$th column.
Let us define the Gram matrix $G$ associated to the
above collection of pseudo-reflexions:
$$G = \left(
\begin{array}{llcl}
t_{00} &t_{10}& \cdots & t_{k-1,0}\\
t_{01} &t_{11}& \cdots & t_{k-1,1}\\
t_{02} &t_{12}& \cdots & t_{k-1,2}\\
\vdots &\vdots& \cdots &\vdots \\
t_{0, k-1} &t_{1, k-1}& \cdots & t_{k-1,k-1}\\
\end{array} \right). \leqno(2.11)$$
We shall treat the cases where $G$ is either symmetric or anti-symmetric.
Let us introduce an upper triangle matrix
$S$ satisfying
$$ G=S+^tS \;\;\;\; (resp.\;G=S-^tS),$$
for a symmetric (anti-symmetric) matrix $G.$
In the anti-symmetric case, we shall use a convention so that the diagonal
part of $S$ is a scalar multiplication on the unit matrix.
It is easy to see that for the symmetric (resp. anti-symemtric)
$G$ the diagonal element $t_{jj}=2$ (resp. $t_{jj}=0$).

\begin {prop}
For an hypergeometric group $\Gamma$ defined over $\bf R,$ the following
statements hold.

1) Suppose that the space of quadratic invariant matrices $X^\Gamma$
is 1 dimensional. Then $X^\Gamma$
coincides with the Gram matrix $G$ (2.11) up
to scalar multiplication.

2)The composition of all generators
$R_0, \cdots, R_{k-1}$ gives us the Seifert form:
$$R_{k-1} \cdots R_0= \mp ^tS \cdot S^{-1}, \leqno(2.12)$$
where to the minus sign corresponds symmetric $G$
and to the plus sign  anti-symmetric $G.$
\label{prop23}
\end{prop}

\noindent
{\bf proof}
1) It is enough to prove that the Gram matrix is a quadratic invariant.
We calculate
$$R_j G ^tR_j = (id_k - Q_j)G(id_k -\;^tQ_j).$$
It is esy to compute
$$ Q_jG= (t_{aj}t_{jb})_{0 \leq a,b \leq k-1 }, G\;^tQ_j
 = (t_{ja}t_{jb})_{0 \leq a,b \leq k-1},$$
$$Q_j G^tQ_j= t_{jj} G\;^tQ.$$
It yields the following equality,
$$ G^tQ_j + Q_j G - Q_jG^tQ_j =
(t_{jb}((1-t_{jj})t_{ja}+t_{aj} ))_{0 \leq a,b \leq k-1 },$$
that  vanishes for $G$ symmetric with $t_{jj}=2$ and for
$G$ anti-symmetric with $t_{jj}=0.$

2)
It is possible to apply directly our situation to that of Theorem
~\ref{thm23}. In the symmetric case, $t_{ii}=2$ and

$$V=\left(
\begin{array}{llcll}
2 & t_{10} & t_{20} &\cdots &t_{k-1,0}\\
0 & 2 &  t_{21} & \cdots &t_{k-1,1} \\
0 &0 & 2 &  \cdots&t_{k-1,2} \\
\vdots &\ddots  & \ddots &\vdots &\vdots\\
0 & 0 & \cdots &  0 &2 \\
\end{array} \right),
U=\left(
\begin{array}{llcll}
0 & 0 & 0 &\cdots &0\\
t_{10} & 0&  0& \cdots &0 \\
t_{20} &t_{21} & 0 &  \cdots&0 \\
\vdots &\ddots  & \ddots &\vdots &\vdots\\
t_{k-1,0} & t_{k-1,1} & \cdots & t_{k-1,k-2}&0 \\
\end{array} \right),$$
in accordance with the notation (2.7).
The formula (2.8) means (2.12) with minus sign.
In the anti-symmetric case $t_{ii}=0, 0 \leq i \leq k-1$
and (2.7) yields (2.11) with plus sign. {\bf Q.E.D.}
\begin{cor}
We can determine the Stokes matrix $S$ by the following relation
$$ S=(id_k -R_{k-1} \cdots R_0)^{-1}G, \leqno(2.13)$$
with the aid of the Gram matrix and pseudo-reflexions.
\label{cor1}
\end{cor}

In some sense, a converse to the Proposition ~\ref{prop23}
holds. To show this, we remember a definition and a proposition from
~\cite{Ok}.
\begin{dfn} The fundamental set $(u_0(\lambda), \cdots u_{k-1}(\lambda))$
of the system $(1.1)"$ is a set of its component solutions
satisfying the following asymptotic expansion:
$$u_j(\lambda)= (\lambda-\lambda_j)^{\rho_j}\sum_{r=0}^\infty g_r^{(j)}
(\lambda-\lambda_j)^r,$$
where $(\lambda_0, \cdots\lambda_{k-1})$ are eigenvalues of the matrix
$\Lambda$ and $\lambda=\lambda_j$ is the unique singular point of
$u_j(\lambda)$ and regular at the  remaining points $\lambda =
\lambda_i,$ $i \not = j$ ($0 \leq j \leq k-1$). The
exponenets $\rho_j$ are diagonal elements of the matrix $A_1(0).$
\label{dfn21} \end{dfn}

After ~\cite{Ok}, the fundamental set to the system $(1.1)"$ is
uniquely determined.

\begin{prop}(~\cite{Ok}) Every generator of an hypergeometric
group $\Gamma$ over $\bf R$ defined for the system of type $(1.1)"$
(without logarithmic solution) is a product of pseudo-reflexions of
the follwing form  expressed with respect to its fundamental set:
$$M_j= id_k -
\left(
\begin{array}{llcllll}
0 & \cdots &0& s_{j0}&0 &\cdots &0\\
0 & \cdots &0& s_{j1}&0 &\cdots &0\\
\vdots &\cdots& \vdots &\vdots &\vdots&\cdots &\vdots\\
0 & \cdots &0& s_{j,k-1}&0 &\cdots &0\\
\end{array} \right).\leqno(2.14)$$
where $s_{jj}=2$ or 0.
\label{prop24}
\end{prop}
We get the following corollary to the above proposition
~\ref{prop24}.

\begin{cor}
Assume that the  hypergeometric group $\Gamma$ is generated by
pseudo-reflexions $T_0, \cdots, T_{k-1}$ such that
$ rank (T_i - id_k)=1$ for $0 \leq i \leq k-1.$ Then it is possible to choose
a suitable set of
pseudo-reflexions generators
$R_j$ like (2.9), (2.10), up to constant multiplication on $Q_j$,
so that they determine the quadratic invariant
Gram matrix
like (2.11).
\label{cor26}
\end{cor}
\noindent
{\bf proof}
The proposition ~\ref{prop24} implies that every generator
$T_i$ is a product of pseudo-reflexions $M_j$ with
$s_{jk}$ possibly different from $t_{jk}.$
From the condition on the quadratic invariant $X^\Gamma$ and the
proposition ~\ref{prop23}, $s_{jk}$ must coincide with $t_{jk}.$
That is to say $\Gamma$ must be generated by
$M_0, \cdots, M_{k-1}$ with $\frac{s_{ja}}{t_{ja}}=\frac{s_{jb}}{t_{jb}} $
for all $a,b,j \in \{0,\cdots, k-1\}.$
This means that $\Gamma$
has as its  generators  the pseudo-reflexions
$R_0, \cdots ,R_{k-1}$
of (2.12) up to constant multiplication on $Q_j$. {\bf Q.E.D.}

\noindent
{\bf proof of Theorem ~\ref{thm12}}
First we remark that solutions to (1.3) have no logarithmic asymptotic
behaviour around any of their singular points except for  the infinity.

In the case with $k$ odd for $X^{H_0},$
there exists $a \not =0$ such that
the vector $\vec{v_0}:=
^t(1+(-1)^{k-1}, -k, \;_{k} C_{2},$ $ \cdots,$
$(-1)^{k-2}\;_{k} C_{k-2},$
$(-1)^{k-1}\;_{k} C_{k-1} )$ $\in {\bf R}^{k}$
satisfies the relation:
$$X^{H_0} \vec{v_0}= ^t(a,0,0,\cdots,0).$$
Actually this fact can be proven almost without calculation in the
follwing way.
First we introduce a series of vetors
$$\vec{w_\ell} = (x_{-\ell}, x_{-\ell+1}, \cdots, x_{k-1-\ell}),
\;\;\;\;
\ell =0,1,\cdots,k-1.$$ Then
the equation $(2.6)''$ can be rewritten in terms of $\vec{w_\ell}$:
$$ \vec{w_\ell} \cdot
\vec{v_0} = \sum_{i=0}^{k-1}(-1)^i\; _kC_i \cdot x_{i-\ell}=0
\;
\rm{for}\;\; 1\leq \ell \leq k-1.$$
On the other hand, the vector
$\vec{w_0}$  is linearly independent of the vectors
$\vec{w_1}, \dots
\vec{w_{k-1}}$ by virtue of the construction of the matrix
$X.$ Therefore $
\vec{w_0} \cdot \vec{v_0} \not= 0$ as  $\vec{v_0} \not= 0.$
This means the existence of the non zero constant $a$ as above.

This relation with  the corollary ~\ref{cor26}
gives immediately the expression below for
the
pseudo-reflexions
$$R_j= id_k -
\left(
\begin{array}{llcllll}
0 & \cdots &0& (-1)^{j+k-1}\;_kC_j \cdot r&0 &\cdots &0\\
\vdots &\cdots& \vdots &\vdots &\vdots&\cdots &\vdots\\
0 & \cdots &0& -(-1)^{k-1}k\cdot r&0 &\cdots &0\\
0 & \cdots &0& (1+(-1)^{k-1})\cdot r&0 &\cdots &0\\
0 & \cdots &0& -k\cdot r&0 &\cdots &0\\
0 & \cdots &0& \;_kC_2 \cdot r&0 &\cdots &0\\
\vdots &\cdots& \vdots &\vdots &\vdots&\cdots &\vdots\\
0 & \cdots &0& (-1)^{k-1-j}\;_kC_{k-j-1}\cdot r &0 &\cdots &0\\
\end{array} \right), \;\; 0\leq j \leq k-1, \leqno(2.15)$$
whose Gram matrix is equal to
$$
\begin{array}{ll}
G_{ij}=(-1)^{i-j+k-1}\;_kC_{i-j}\cdot r\;\;\;\;& i >j\\

G_{ii}= (1+(-1)^{k-1})\cdot r\;\;\;\; &i = j\\

G_{ij}=(-1)^{i-j}\;_kC_{j-i}\cdot r\;\;\;\;& j >i\\
\end{array}
\leqno(2.16)$$
with some constant $r$.
As for the case $k-$ even, the equations $(2.6)'''$ and the Corollary ~\ref{cor26}
gives us the expression (2.15) for the pseudo-reflexion generators.

Taking into account the theorem ~\ref{thm11} for the symmetric Gram matrix, we obtain the desired
statement for the case $k-$ odd, as it is required from
Proposition ~\ref{prop24} $G_{ii}=2=2r.$

For the case $k-$even,
we remember a statement on the Stokes matrix
from \cite{BJL}( Proposition 1,2) 
which claims that if the matrix $\mu$ of (1.1) has
integer eigenvalues, the equality $ det(S + \;^tS)=0$ must hold.
The Corllary~\ref{cor1} gives us the relation
$$S= (id_k + (id_k-V)(id_k+U)^{-1})^{-1}G= (id_k+U)G^{-1}G=id_k+U,$$
with
$$ U_{ij}=(-1)^{i-j+k-1}\;_kC_{i-j}\cdot r\;\;\;\; i >j.$$
We shall choose the constant $r=1$ so that $S+ \;^tS= 2 id_k + U + ^tU$
posesses an eigenvector $(1,-1,\cdots, 1,-1)$ with zero eigenvalue.
{\bf Q.E.D.}
\begin{remark}
{\em The Gram matrix (2.16) that has been calculated for the
fundamental set (Definition ~\ref{dfn21}) of the equation (1.3) gives
directly a suitable Stokes matrix we expected. For other Fano
varieties, however, the Gram matrix calculated with respect to the
fundamental set does not necessarily give a desirable form, as it is
seen from the case of odd dimensional quadrics.  This situation makes
us to be careful in the choice of the base of solutions for which we
calculate the Gram matrix.} \end{remark}
{ \center{\section{\bf
Geometric interpretation of the hypergeometric equation
}}
}

In this section we show that the equation (1.4) arises from
the differential operator that annihilates the fibre integral
associated to the family of variety defined as a complete
intersection
$$ X_s := \{(x_0, \cdots x_k) \in {\bf C}^{k+1};
f_1(x) + s=0, f_2(x)+1=0\}.
\leqno(3.1)$$
where
$$f_1(x)= x_0 x_1 \cdots x_k, f_2(x)=x_0+ x_1 + \cdots +x_k.$$
This result has been already
announced by \cite{Giv}, \cite{Gol} and \cite{Bar}.
Our main theorem of this section is the following
\begin{thm}
Let us  assume that $\Re (f_1(x) +s) |_{\Gamma} <0,$
$\Re (f_2(x)+s) |_{\Gamma } <0,$ out of a compact set
for a Leray coboundary
cycle $\Gamma \in H^{k+1}({\bf C}^{k+1}\setminus X_s)$
avoiding the hypersurfacess $f_1(x) + s=0$ and $f_2(x)+1=0.$
For such a cycle we consider the following residue integral:
$$ I^{(v_1, v_2)}_{x^\bi, \Gamma}(s)=
\int_\Gamma x^{\bi+\b1}(f_1(x) + s)^{-v_1}
(f_2(x) +1)^{-v_2}\frac{dx}{x^{\b1}},
\leqno(3.2)
$$ for the monomial 
$x^\bi := x^{i_0}_0\cdots x^{i_{k}}_k, 
x^\b1 :=
x_0\cdots x_k.$ 
Then the  integral $I^{(1, 1)}_{x^{0},\Gamma}(s)$
satisfies the following hypergeometric
differential equation
$$ [\vartheta_s^k - k^ks(\vartheta_s + \frac{1}{k})
( \vartheta_s + \frac{2}{k})\cdots
(\vartheta_s + \frac{k}{k})]I^{(1, 1)}_{1, \Gamma}(s)=0 \leqno(3.3)$$
which has unique holomorphic solution at $s=0,$
$$ I_0(s)= \sum_{m\geq 0} \frac{(km)!}{(m!)^k} s^m.   \leqno(3.4)$$
\label{thm31}
\end{thm}
We shall put $\zeta = \frac{1}{k^k s},$ to get (1.4) from (3.3).
Our calculus is essentially based on the Cayley trick method developed
in \cite{T}.

{\bf Proof of Theorem 3.1}
Let us consider the Mellin transform of the fibre integral (3.2)
$$ M_{{\bi},\Gamma}^{(v_1, v_2)} (z):=\int_\Pi s^z
I^{(v_1, v_2)}_{x^\bi}(s)
\frac{ds}{s}
\leqno(3.5)$$
For the Mellin transform (3.5), we have the following
$$ M_{{\bi},\Gamma}^{(v_1, v_2)} (z)= g(z)
\prod^{k-1}_{\ell =0}\Gamma(z+ i_\ell+1
-v_2) \Gamma( - \sum_{\ell =0}^{k-1}(i_\ell+1) -kz + v_1 + kv_2)\Gamma(-z+ v_2)
\Gamma(z), \leqno(3.5)'$$
with  $g(z)$ a rational function in $e^{\pi i z}.$
The formula $(3.5)'$ shall be proven below.
In substituting $\bi = 0, v_1=v_2=1,$
we see that
$$I^{(1, 1)}_{x^{0}, \Gamma}(s)= \int_{\check\Pi} s^{-z}g(z)
\frac{\Gamma(z)^k} {\Gamma(kz)} dz,$$
where $\check\Pi$ denotes the path $(-i\infty, +i\infty)$ avoidoing the poles of $\Gamma(z)= 0, -1, -2, \cdots.$
From this integral representation, the equation
$(3.3)$ immediately follows in taking account the fact that
the factor $g(z)$ plays no role in establishment of the
differential equation. {\bf Q.E.D.}

{\bf Proof of $(3.5)'$}
In making use of the Cayley trick, we transform the integral
$(3.5)$ into the following form
$$M_{{\bi},\Gamma}^{(v_1, v_2)} (z)= 
\int_{\Pi \times \bar {\bR }_+^2 \times \Gamma}x^{\bi+\b1} 
e^{y_1(f_1(x) +s)+ y_2(f_2(x)+1)} y_1^{v_1} y_2^{v_2} s^z 
\frac{dx}{x^\b1}\frac{dy}{y^\b1}
\frac{ds}{s^\b1}, \leqno(3.6)$$
with $\bar \bR_+$ a double covering cycle of the positive real axis 
in $\bC_{y_p}$  turning around $y_p=0$ 
that begins and returns to $\Re y_p \rightarrow \infty,$ for $p=1$ or $2.$
Here we introduce new variables $T_0, \cdots T_{k+2},$
$$ T_i = y_1 x_i, \;\; 0 \leq i \leq k-1, T_k =y_1s, T_{k+1}= y_2
x_0x_1 \cdots x_{k-1}, T_{k+2}=y_2 \leqno(3.7)$$
in such a way that the phase function of the right hand side
of (3.6) becomes
$$ y_1(f_1(x) +s)+ y_2(f_2(x)+1)=T_0+ T_1 + \cdots + T_{k+2}.$$
If we set
$$Log\; T := ^t(log\; T_0, \cdots, log\; T_{k+2})$$
$$ \Xi := ^t(x_0, \cdots, x_{k-1},  s, y_1,y_2) $$
$$Log\; \Xi := ^t(\log\; x_0, \cdots, \log\; x_{k-1},  \log\; s,
,\log\; y_1, log\; y_2).$$
Then the above relationship (3.7) can be written down as
$$ Log\; T= {\sf L}\cdot Log\; \Xi, \leqno(3.8)$$
where
$${\sf L}= \left [\begin {array}{cccccccc}
\noalign{\medskip}1&0&0&\cdots&0&0&1&0\\
\noalign{\medskip}0&1&0&\cdots &0&0&1&0\\
\noalign{\medskip}0&0&1&\cdots&0&0&1&0\\
\noalign{\medskip}\vdots&\vdots&\vdots&\ddots&\vdots&\vdots&\vdots&\vdots\\
\noalign{\medskip}0&0&0&\cdots&1&0&1&0\\
\noalign{\medskip}0&0&0&\cdots&0&0&1&0\\
\noalign {\medskip} 1&1&1&\cdots&1&0&0&1\\
\noalign{\medskip}0&0&0&\cdots&0&1&0&1
\end {array}\right ].$$
This yields immediately
$$Log\; \Xi= {\sf L}^{-1}\cdot Log\; T,$$
with
$${\sf L}^{-1}= \left [\begin {array}{cccccccc}
1&0&0&\cdots&0&-1&0&0\\
\noalign{\medskip}0&1&0&\cdots &0&-1&0&0\\
\noalign{\medskip}0&0&1&\cdots&0&-1&0&0\\
\noalign{\medskip}\vdots&\vdots&\vdots&\ddots&\vdots&\vdots&\vdots&\vdots\\
\noalign{\medskip}0&0&0&\cdots&1&-1&0&0\\
\noalign{\medskip}1&1&1&\cdots&1&-k&-1&1\\
\noalign {\medskip} 0&0&0&\cdots&0&1&0&0\\
\noalign{\medskip}-1&-1&-1&\cdots&-1&k&1&0
\end {array}\right ].$$
If we set
$$ (i_0, \cdots i_{k-1},z , v_1, v_2)\cdot{\sf L}^{-1}
= \bigl ({\mathcal L}_0(\bi, z, v_1, v_2), \cdots,
 {\mathcal L}_{k+2}(\bi, z, v_1, v_2) \bigr). \leqno(3.9)$$
then we can see that
$$M_{{\bi},\Gamma}^{(v_1, v_2)} (z)=
\int_{\Pi  \times
{\bar \bR}_+^2 \times \Gamma}x^{\bi+\b1} e^{T_0 + \cdots +T_{k+2}}
y_1^{v_1} y_2^{v_2} s^z
\frac{dx}{x^\b1}\frac{dy}{y^\b1}
\frac{ds}{s^\b1}$$
$$= \int_{\sf L_{\ast}(\Pi  \times {\bar\bR}_+^2 \times \Gamma)} e^{T_0
+ \cdots +T_{k+2}} \prod_{0\leq i \leq k+2} T_i^{{\mathcal L}_{i}(\bi,
z, v_1, v_2)}\bigwedge_{0\leq i \leq k+2} \frac {dT_i}{T_i}.$$ Here
${\sf L}_{\ast}(\Pi  \times
{\bar\bR}_+^2 \times \Gamma)$ denotes a
$(k+3)-$cycle in  
$T_0  \cdots T_{k+2}\not =0$ that obtained as a
image of
$\Pi  \times {\bar\bR}_+^2 \times \Gamma$ under the
transformation induced by
$\sf L.$ 
In view of the choice of the cycle
$\Gamma,$ we can apply the formula to calculate $\Gamma$ function to
our situation:  $$ \int_{C} e^{-T}T^\sigma \frac{dT}{T} = (1- e^{2 \pi
i \sigma}) \Gamma(\sigma),$$ for the unique nontrivial  cycle $C$
turning around $T=0$ that begins and returns to $\Re T \rightarrow +
\infty.$ Here one can consider the natural action $\lambda:C_a
\rightarrow \lambda(C_a)$ defined by the relation,
$$\int_{\lambda(C_a)} e^{-T_a}T_a^{\sigma_a} \frac{dT_a}{T_a} =
\int_{(C_a)} e^{-T_a}(e^{2\pi \sqrt -1 }T_a)^{\sigma_a}
\frac{dT_a}{T_a}.$$ In terms of this action ${\sf L}_\ast (\Pi  \times
{\bar\bR}_+^2 \times \Gamma) $ is shown to be homologous to a chain
$$\sum_{(j_1^{(\rho)},\cdots,j_L^{(\rho)} )\in [1,\Delta]^L} 
m_{j_1^{(\rho)},
\cdots,j_L^{(\rho)}} \prod_{a=1}^L \lambda^{j_a^{(\rho)}}(C_a),$$ with
$m_{j_1^{(\rho)}, \cdots, j_L^{(\rho)}}  \in \bZ. $  
This explains the
appearance of the factor $g(z)$ in front of the $\Gamma$ function
factors in $(3.5)'$.

The direct calculation of (3.9) shows that
$$ {\mathcal L}_{\ell}(\bi, z, v_1, v_2)= z+ i_\ell+1-v_2,
0 \leq \ell \leq k-1,$$
$${\mathcal L}_{k}(\bi, z, v_1, v_2)= - \sum_{\ell =0}^{k-1}(i_\ell+1)+
v_1 + k(v_2-z), {\mathcal L}_{k+1}(\bi, z, v_1, v_2)=-z+ v_2,
{\mathcal L}_{k+2}(\bi, z, v_1, v_2)=z.$$
This shows the formula $(3.5)'$. {\bf Q.E.D.}

 In combining Theorems ~\ref{thm12}, ~\ref{thm31}, we can state that we
found out a deformation of an algebraic variety $X_\lambda=\{
(\frac{\lambda}{k})^k(x_0 x_1 \cdots x_k) + 1=0, 
x_0+ x_1 + \cdots +x_k=1=0\}.$
such that its variation gives rise to the equation (1.3).
It means that we establish a connexion between an exceptional collection
of ${\bf CP}^{k-1}$ and a set of vanishing cycles for its mirror
counter part $X_\lambda.$  Thus our theorems give an affirmative answer
to the hypothesis stating the existence of
such relationship between two mirror symmetric varieties (so
called Bondal-Kontsevich  hypothesis) in a special case.
See ~\cite{Gol} and ~\cite{Konts} in this
respect for the detail.

It is known from the theory  of period integrals associated to the
complete intersections ~\cite{Gre} that the integrals
$I^{(v_1, v_2)}_{x^{\bi}, \Gamma}((\frac{k}{\lambda})^k)$
for $\Gamma_\lambda \in H_{k+1}(\bC^{k+1} \setminus X_\lambda, \bZ)$
such that the cycle
$\gamma_\lambda \in H_{k-1}(
X_\lambda, \bZ)$ (of which Leray's coboundary is 
$\Gamma_\lambda$)
has singularities only at the discriminant locus of $X_\lambda$
where the cycle $\gamma_\lambda $ becomes singular (or vanishes).
On the other hand, in \S 2 we found a set of solutions
such that $u_j(\lambda)$ has an unique singular point
$\lambda=e^{2\pi \sqrt -1 \frac{j}{k}}$ and regular at the  remaining
points $\lambda = e^{2\pi \sqrt -1 \frac{i}{k}},$ $i \not = j.$  Two
solutions to an hypergeometric differential equation $(1.3)$ with the
same assigned asymptotic behaviours at all possible singular points
must coincide. In combination of this argument with the
Picard-Lefschetz theorem, we obtain the following.

\begin{cor}
There exists a set of cycles  $\gamma_j \in H_{k-1}(
X_\lambda, \bZ),$ $0 \leq j \leq k-1$ such that for their Leray's
coboundary $\Gamma_j \in H_{k+1}(\bC^{k+1} \setminus X_\lambda, \bZ)$
we have the identity,
$$I^{(1, 1)}_{x^{0}, \Gamma_j}((\frac{k}{\lambda})^k)=
u_j(\lambda), \;\;\; 0 \leq j \leq k-1,
$$
with $u_j(\lambda)$ the fundamental solution to $(1.3)$ in the sense of
Definition ~\ref{dfn21}. Consequently the Gram matrix $G$ of $(2.16)$
is equal to the intersection matrix $(<\gamma_i,\gamma_j>)_{0 \leq i,j
\leq k-1}$ after proper choice of constant $r=1.$
\label{cor32}
\end{cor}


\vspace{\fill}


%

\noindent

\begin{flushleft}
 \begin{minipage}[t]{6.2cm}
  \begin{center}
{\footnotesize Independent University of Moscow\\
Bol'shoj Vlasijevskij Pereulok 11,\\
  MOSCOW, 121002,\\
Russia\\
{\it E-mails}:  tanabe@mccme.ru,}
\end{center}
\end{minipage}\hfill
\end{flushleft}
\end{document}